\documentclass[11pt]{article}
\usepackage{graphicx}
\usepackage[letterpaper]{geometry}
\usepackage{amsmath}
\usepackage{latexsym}

\newtheorem{theorem}{Theorem}
\newtheorem{ob}[theorem]{Observation}
\newtheorem{lemma}[theorem]{Lemma}
\newtheorem{corollary}[theorem]{Corollary}
\newtheorem{conjecture}{Conjecture}

\newenvironment{Proof}{\proofStarter{Proof.}}{\hspace{8mm}\QEDmark\smallskip}
\newcommand{\QEDmark}{\mbox{\Large $\diamond$}}
\newcommand{\proofStarter}[1]{\textsc{#1} }

\newcommand{\BREAK}[2]{\mathop{\mathit{b}}_{#1}(#2)}

\newcommand{\cor}{\mathop{\mathit{cor}}}

\begin{document}

\begin{center}
{\large
\textbf{Vertex Colorings without Rainbow or Monochromatic Subgraphs} \bigskip \\ }
Wayne Goddard and Honghai Xu \smallskip \\
Dept of Mathematical Sciences, Clemson University \\
Clemson SC 29634 \\
\textsf{\{goddard,honghax\}@clemson.edu}
\end{center}

\begin{abstract}
This paper investigates vertex colorings of graphs such that 
some rainbow subgraph~$R$ and some monochromatic subgraph $M$ are
forbidden. Previous work focussed on the case that $R=M$. Here we
consider the more general case, especially the case that $M=K_2$.
\end{abstract}

\section{Introduction}

Let $F$ be a graph. Consider a coloring of the vertices of $G$.
We say that a copy of $F$ (as a subgraph) is \emph{rainbow} (or heterochromatic)
if all its vertices receive different colors.
We say that the copy of $F$ is \emph{monochromatic}
if all its vertices receive the same color.

The question of avoiding monochromatic copies of a graph is well
studied (see for example the survey \cite{Tuza97local}). Less studied, but still common,
is the question of avoiding rainbow copies (especially for edge-colorings); 
see for example \cite{AI08rainbow,BSTDP12stars,BSTSD10consecutive}.
In \cite{GWX15firstWorm,GWX14secondWorm} we defined WORM colorings: these forbid 
both a rainbow and a monochromatic
copy of a specific subgraph. But it is more flexible to 
allow different restrictions. For graphs $M$ and $R$, we define an 
\emph{$(M, R)$-WORM coloring of $G$} to be a coloring of the vertices of 
$G$ with neither a monochromatic subgraph 
isomorphic to $M$ nor a rainbow subgraph isomorphic to $R$. 
Note that such a coloring is not guaranteed to exist. For example, any $G$ with at least one edge
does not have a $(K_2,K_2)$-WORM coloring. 

This coloring is a special case of the
``mixed hypergraphs'' introduced by Voloshin (see for example~\cite{voloshin1995upper});
see~\cite{TV08problems} for an overview.
A related question studied in the edge case is 
the \emph{rainbow Ramsey number} (or \emph{constrained Ramsey number});
this is 
defined as the minimum~$N$ such that any coloring of the edges of $K_N$
produces either a monochromatic~$M$ or a rainbow $R$. See~\cite{Eroh04ramsey}. 

One special case of WORM colorings has a distinguished history. Erd\H{o}s et al.~\cite{EHKRS86locally}
defined the \emph{local chromatic number} of a graph as the maximum
order of a rainbow star that must appear in all proper colorings.
In our notation, this is the minimum $r$ such that the graph has an $(K_2, K_{1,r+1})$-WORM
coloring. For a survey on this parameter, see \cite{Osang13local}.

One case is trivial: if we forbid a rainbow $K_2$, then every component
of the graph must be monochromatic. Similarly, if we forbid a rainbow $kK_1$,
then this is equivalent to using less than $k$ colors. 
So we will assume that the subgraph $R$ has
at least three vertices and at least one edge. On the other hand, taking $M=K_2$
is equivalent to insisting that the coloring is proper.
Also, taking $M=kK_1$ is equivalent to using each color less than $k$ times.

Having two competing restrictions leads naturally to considering both the 
minimum and maximum number of colors in such a coloring. So we define
the upper chromatic number
$W^+(G; M, R)$ as the maximum number of colors,  
and the lower chromatic number $W^-(G; M, R)$ as the minimum number of colors, in an $(M, R)$-WORM 
coloring of~$G$ (if the graph has such a coloring). For bounds, it will be useful
to also let $m^-(G;M)$ be the minimum number of colors without a monochromatic $M$,
and $r^+(G;R)$ be the maximum number of colors without a rainbow $R$. Note that
\[
   m^-(G;M) \le W^-(G; M, R) \le W^+(G; M, R) \le r^+(G;R) ,
\]
provided $G$ has an $(M,R)$-WORM coloring.

We proceed as follows.
Section~\ref{s:prelim} contains some general observations.
In Section~\ref{s:path} we provide one general upper bound
when $R$ is a path. In Section~\ref{s:proper} we consider
proper colorings without rainbow $P_3$, $P_4$, or $C_4$. 
Finally, in Section~\ref{s:other} we provide a few results for other cases.

\section{Preliminaries} \label{s:prelim}

We start with
some simple observations.
If $G$ is bipartite then the bipartition is immediately an 
$(M,R)$-WORM coloring.
Indeed, if $G$ is $k$-colorable with $k<|R|$ , then a 
proper $k$-coloring of $G$ is an $(M,R)$-WORM coloring. Also:

\begin{ob} Fix graphs $M$ and $R$ and let $G$ be a graph.

(a) If $G$ has an $(M,R)$-WORM coloring, then 
so does $G-e$ where $e$ is any edge and $G-v$ where $v$ is any vertex.
Further, $W^+(G-e; M,R) \ge W^+(G; M,R)$
and $W^+(G-v; M,R) \ge W^+(G; M,R)-1$, with similar results for the
lower chromatic number.

(b) 
If $M$ and $R$ are connected but
$G$ is disconnected, then $W^+(G; M,R)$ is the sum of the parameter for
the components, and $W^-(G; M,R)$ is the maximum of the parameter
for the components.

(c) It holds that 
$W^+(G; M,R) = |V(G)|$ if and only if $G$ is $R$-free.

(d) It holds that 
$W^+{(G;M,R}) \ge |R|-1$ if $G$ is $|R|-1$ colorable (and has at least that many vertices).
\end{ob}

We will also need the following idea from~\cite{GX15rainbow}.
We say that a set $S$ \emph{bi-covers}
a subgraph~$H$ if at least two vertices of $H$ are in $S$. For positive integer $s$,
define $\BREAK{F}{s}$ to be the maximum number of copies of $F$ that can be bi-covered by
using a set of size $s$. (Note that by definition $\BREAK{F}{1}=0$.) 

\begin{lemma} \label{l:breakBound} \cite{GX15rainbow}
Suppose that graph $G$ of order $n$ contains $f$ copies of $R$ and that
$\BREAK{R}{s} \le a(s-1)$ for all $s$. Then $r^+(G;R) \le n - f/a$.
\end{lemma}

\subsection{General $M$}

It should be noted that maximizing the number of colors while avoiding a rainbow
subgraph can produce a large monochromatic subgraph. For example:

\begin{ob}
For all connected graphs $M$, there exists a graph $G$ such that 
$W^+(G; M, P_3) < r^+(G; P_3)$.
\end{ob}
\begin{Proof}
In~\cite{GX15rainbow} we considered the 
\emph{corona} $\cor(G)$ of a graph $G$; this is the graph obtained from $G$ by adding,
for each vertex $v$ in $G$, a new vertex $v'$ and the edge $vv'$. 
It was shown that $r^+(G; P_3) = |G|+1$. In fact, we note here that if $G$ is connected,
then one can 
readily show by induction that the optimal coloring is unique and gives every vertex of $G$ 
the same color. In particular, it follows that the no-rainbow-$P_3$ coloring
of $\cor(M)$ with the maximum number of colors contains a monochromatic copy of $M$. 
\end{Proof}

\section{A Result on Rainbow Paths} \label{s:path}

We showed~\cite{GWX15firstWorm} that a nontrivial graph $G$ has a $(P_3,P_3)$-WORM coloring
if and only if it has a $(P_3,P_3)$-WORM coloring using only two colors. 
We prove an analogue for general paths. This result is a slight generalization of Theorem~10 in~\cite{TV00uncolorable}.

\begin{theorem} \label{t:pathUpper} Fix some graph $M$; 
if graph $G$ has an $(M,P_r)$-WORM coloring, then $G$ has one using at most $r-1$ colors.
\end{theorem}
\begin{Proof}
Consider an $(M,P_r)$-WORM coloring $f$ of $G$. Let $G_M$ be the spanning subgraph of $G$
whose edges are monochromatic and $G_R$ the spanning subgraph whose
edges are rainbow. It follows that $G_M$ that does not contain $M$,
and that $G_R$ does not contain $P_r$. It is well known that 
a graph without $P_r$ 
has chromatic number at most $r-1$. 
 
Now, let $g$ be a proper coloring of $G_R$ using at most $r-1$ colors and consider
$g$ as a coloring of $G$. Note that 
the monochromatic edges under $g$ are a subset of those under $f$.
Therefore, $g$ is a $(M,P_r)$-WORM coloring of $G$ using at most $r-1$ colors.
\end{Proof}

It follows that:

\begin{corollary} \label{cor:pathW-}
For any graph $M$ and $r>0$, graph
$G$ has an $(M,P_r)$-WORM coloring if and only if $m^-(G,M) \le r-1$.
If so, $W^-(G; M, P_r ) = m^-(G,M)$.
\end{corollary}

On the other hand, Theorem~\ref{t:pathUpper} does not extend to stars. For example,
Erd\H{o}s et al.~\cite{EHKRS86locally} constructed a shift graph that 
has arbitrarily large chromatic number
but can be properly colored without a rainbow $K_{1,3}$. That is:

\begin{theorem} \label{t:erdos} For $r\ge 3$ and $k\ge 1$,
there is a graph $G$ with $W^- ( G; K_2, K_{1,r} ) \ge k$.
\end{theorem}

Nor does Theorem~\ref{t:pathUpper} generalize to $K_3$; see \cite{XuTHesis} and \cite{BTk3worm}.

\section{Proper Colorings} \label{s:proper}

Recall that $W^+(G; K_2,R)$
is the maximum number and $W^-(G; K_2,R)$ the minimum number of colors in 
a proper coloring without a rainbow $R$.

\subsection{Two simple cases}
 
Two cases for $M=K_2$ are immediate:

\begin{ob}
A graph $G$ has a $(K_2,P_3)$-WORM  coloring if and only if it is bipartite.  
If so, $W^+(G: K_2,P_3) = W^-(G; K_2,P_3)=2$, provided $G$ is connected and nonempty.
\end{ob}
\begin{Proof}
If we have a $(K_2,P_3)$-WORM coloring, then for each vertex $v$ 
all its neighbors must have the same color, which is different to $v$'s color. 
It follows that every path must alternate colors.
\end{Proof}

In a proper coloring of a graph, all cliques are rainbow. Thus it follows:

\begin{ob}
A graph $G$ has a $(K_2,K_m)$-WORM coloring if and only if it is $K_m$-free.
If so, $W^+(K_2,K_m) = |G|$ while $W^-(K_2,K_m)$  is the chromatic number of $G$.
\end{ob}

\subsection{No rainbow $K_{1,3}$}

Consider first that $G$ is bipartite. Then in maximizing the colors,
it is easy to see that one may assume the colors in the 
partite sets are disjoint. (If red is used in both partite sets,
then change it to pink in one of the sets.) In particular,
unless $G$ is a star, one can use at least two colors in each 
partite set. (This result generalizes to $R$ any star.)
For example, it follows that $W^+( K_{m,m}; K_2, K_{1,3} ) = 4$
for $m\ge 2$.

Indeed, it is natural to consider the \emph{open neighborhood hypergraph}
$ON(G)$ of the graph~$G$. This is the hypergraph with vertex set $V(G)$ and a hyperedge
for every open neighborhood in $G$. In general, since we have a proper coloring, the 
requirement of no rainbow $K_{1,r}$ is equivalent to every
hyperedge in $ON(G)$ receiving at most $r-1$ colors.
In the case that $G$ is bipartite, the two problems are equivalent:

\begin{ob} \label{o:onh}
For any graph $G$, the parameter $W^+(G; K_2, K_{1,r})$ is at most the 
maximum number of colors in a coloring of $ON(G)$ with
every hyperedge receiving at most $r-1$ colors. Furthermore,
there is equality if $G$ is bipartite.
\end{ob}
\begin{Proof}
When $G$ is bipartite, the $ON(G)$ can be partitioned into two disjoint hypergraphs
and so will have disjoint colors in the hypergraphs. It follows that 
that coloring back in $G$ will be proper.
\end{Proof}

Recall that a $2$-tree is defined by starting with $K_2$ and repeatedly adding
a vertex that has two adjacent neighbors. For example, this includes
maximal outerplanar graphs.

\begin{ob}
If $G$ is a $2$-tree of order at least $3$, then $W^+(G; K_2, K_{1,3})=3$.
\end{ob}
\begin{Proof}
Any $2$-tree is $3$-colorable. Furthermore, it follows readily by 
induction that  a $(K_2, K_{1,3})$-WORM coloring
can use only three colors: when we add a vertex $v$ 
and join it to adjacent vertices $x$ and $y$, they already
 have a common neighbor $z$,
and so $v$ must get the same color as $z$.
\end{Proof}

Osang showed that determining whether a graph has a $(K_2,K_{1,3})$-WORM
coloring
is hard:

\begin{theorem} \cite{Osang13local} 
Determining whether a graph has a $(K_2,K_{1,3})$-WORM coloring
is NP-complete.
\end{theorem}

\subsubsection{Cubic graphs}

We consider next $3$-regular graphs.
Since cubic graphs (other than $K_4$) are $3$-colorable, they have a
$(K_2,K_{1,3})$-WORM coloring. And that coloring uses at most three colors.
Further, they have a coloring using two colors if and only if they are bipartite.
So the only interesting question is the behavior of the upper chromatic number.

\begin{ob} \label{o:cubicClaw}
If $G$ is cubic of order $n$, then $W^+( G; K_2, K_{1,3} ) \le 2n/3$.
\end{ob}
\begin{Proof}
Since $G$ is cubic, the hypergraph $ON(G)$ is $3$-regular and $3$-uniform.
Further we need a coloring of $ON(G)$ where every hyperedge
has at least one pair of vertices the same color.
Consider some color used more than once, say red. 
If there are $r$ red vertices, then at most $3r/2$ hyperedges can have two
red vertices. (Each red can be used at most thrice.) 

It follows that if the $i$\textsuperscript{th} non-unique color is used $r_i$ times, then
we need $\sum_i r_i \ge 2n/3$. Let $B$ be the number of vertices that can 
be discarded and still have one vertex of each color. Then $B=\sum_i (r_i-1)$
and by above $B \ge n/3$. It follows that the total number of colors is at most $2n/3$.
\end{Proof}

\begin{figure}[!ht]
\centerline{\includegraphics{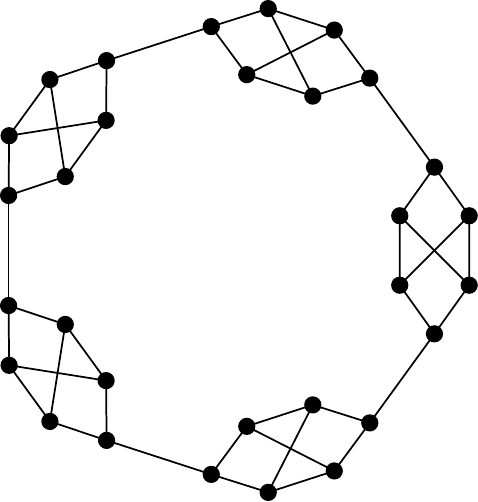}}
\caption{A cubic graph $G$ with $W^+(G;K_2,K_{1,3})$ two-thirds its order}
\label{f:bracelet}
\end{figure}

Equality in Observation~\ref{o:cubicClaw} is obtained by taking disjoint copies of 
$K_{3,3}-e$ and adding
edges to make the graph connected. See Figure~\ref{f:bracelet}.

\begin{figure}[!ht]
\centerline{
\begin{tabular}{c}{\includegraphics{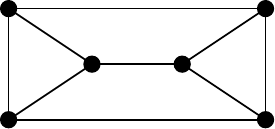}}\end{tabular}
\quad
\begin{tabular}{c}{\includegraphics{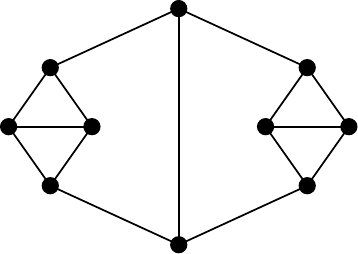}}\end{tabular}
\quad
\begin{tabular}{c}{\includegraphics{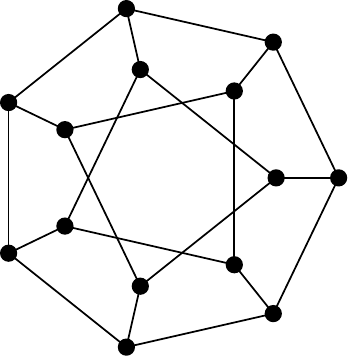}}\end{tabular}}
\caption{The known cubic graphs with $W^+(G;K_2,K_{1,3}) = 3$}
\label{f:cubicThree}
\end{figure}

Consider next the minimum value of $W^+(G; K_2, K_{1,3})$ for cubic graphs of order~$n$. We noted above 
that 
bipartite graphs in general have a value of at least $4$.
Computer search shows that this parameter is at least $3$ for $n\le 18$. Indeed,
it finds only three graphs where the parameter is $3$: one of order $6$ (the prism),
one of order $10$, and one of order $14$,  
the generalized
Petersen graph. These three graphs are shown in Figure~\ref{f:cubicThree}.

It is unclear what happens in general.

\subsection{Forbidding rainbow $P_4$}

We consider proper colorings without rainbow $P_4$'s.
Theorem~\ref{t:pathUpper} applies. That is, a graph $G$ has 
a $(K_2,P_4)$-WORM coloring
if and only if $G$ has chromatic number at most $3$. In particular,
this means that it is NP-complete to determine if a graph has a 
$(K_2,P_4)$-WORM coloring.
Further, if such a coloring exists, then  $W^-(G; K_2, P_4)$ is the 
(ordinary) chromatic number of $G$.
So we consider only the upper chromatic number here.

\begin{ob} \label{o:bipartiteP4}
If graph $G$ is bipartite of order $n$, then $W^+(G; K_2, P_4) \ge n/2+1$.
\end{ob}
\begin{Proof}
In the smaller partite set, give all vertices the same color,
and in the other partite set, give all vertices unique colors. Note that
every copy of $P_4$ contains two vertices from both partite sets.
\end{Proof}

\begin{ob} \label{o:matchingP4}
If connected graph $G$ of order $n$ has a perfect matching, then it holds that
$W^+(G; K_2, P_4) \le n/2+1$.
\end{ob}
\begin{Proof}
Number the edges of the perfect matching $e_1, \ldots, e_{n/2}$ such that for all $i>1$, at least
one of the ends of $e_i$ is connected to some $e_j$ for $j<i$. Then $e_i$, $e_j$, and the 
connecting edge form a $P_4$. It follows that $e_j$ and $e_i$ share a color. Thus the total
number of colors used is at most $2+(n/2-1) = n/2+1$.
\end{Proof}

\begin{figure}[!ht]
\centerline{\includegraphics{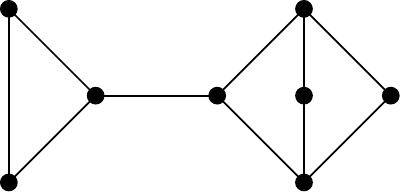}}
\caption{A nonbipartite graph $G$ with a perfect matching and maximum $W^+(G;K_2,P_4)$}
\label{f:matchingEquality}
\end{figure}

For example, equality is obtained in both observations for any connected bipartite graph 
with a perfect matching, such as the balanced complete  bipartite graph or the path/cycle of even
order. Equality is also obtained in Observation~\ref{o:bipartiteP4} for the tree of diameter three 
where the two central vertices have the same degree. Also, there are nonbipartite graphs that 
achieve equality in Observation~\ref{o:matchingP4}; for example, the graph shown in 
Figure~\ref{f:matchingEquality}.

We determine next the parameter for the odd cycle:

\begin{ob}  If $n$ is odd, then 
$W^+(C_n; K_2, P_4)$ is $3$ for $n\le 5$,
and $(n-1)/2$ for $n\ge 7$.
\end{ob}
\begin{Proof}
The result for $n=3$ is trivial and for $n=5$ is easily checked. 
So assume $n\ge 7$. For the lower bound, color red a maximum
independent set, give a new color to every vertex with two red neighbors,
and color each vertex with one red neighbor the same color as on the other
side of its red neighbor. For example, the coloring for $C_{13}$ is shown in 
Figure~\ref{f:cycleColoring} (where the red vertices are shaded).

\begin{figure}[!ht]
\centerline{\includegraphics{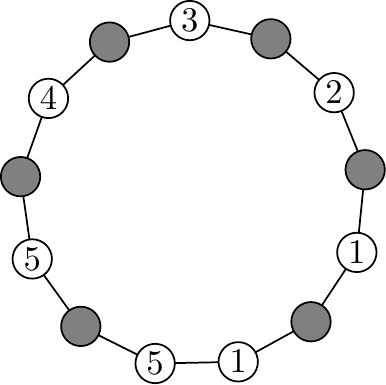}}
\caption{Coloring showing $W^+(C_{13};K_2,P_4)$}
\label{f:cycleColoring}
\end{figure}

We now prove the upper bound.
Two same-colored vertices distance $2$ apart bi-cover two copies of~$P_4$,
while two same-covered vertices distance $3$ apart bi-cover one copy. It follows that if a color is used
$k$ times, it can bi-cover at most $2(k-1)$ copies of $P_4$, except if the vertices
of that color form a maximum 
independent set,
when it bi-covers $2k-1$ copies. Since there are $n$ copies of $P_4$ in total,
by Lemma~\ref{l:breakBound} it follows that the total number of colors is at most
$n/2$, unless some color is a maximum independent set. So say red is a maximum independent
set. Let $b$ and $e$ be the two red vertices at distance~$3$; say the portion of the
cycle containing them is $\mathit{abcdef}$.
By considering the $a$--$d$ copy of $P_4$,
it follows that $a$ must have the same color as $c$ or~$d$. Similarly, $f$ must
have the same color as $c$ or~$d$. Thus the total number of colors is at most
$1+(n-(n-1)/2)-2 = (n-1)/2$. 
\end{Proof}

In contrast to Observation~\ref{o:bipartiteP4}, we get the following:

\begin{theorem}  \label{t:pathTriangle}
If connected graph $G$ has every vertex in a triangle, then 
 $W^+(G; K_2,P_4)= 3$ if such a coloring exists.
\end{theorem}
\begin{Proof}
Note that every triangle is properly colored.
We show that every triangle receives the same three colors.
Consider two triangles $T_1$ and $T_2$. If $T_1$ and $T_2$ 
share two vertices, then the third vertex in 
each share a color. Consider the case that $T_1$ and~$T_2$
share one vertex. Then by considering the four $P_4$'s using all vertices
but one, it readily follows that the triangles must have
the same colors.

Now, assume that $T_1$ and $T_2$ are disjoint but joined by an edge $e$.
 Suppose they do not have the
same three colors. Then there is vertex $u_1$ in $T_1$ and $u_2$ in $T_2$
that do not share a color with the other triangle. If $u_1$ and $u_2$ are the ends
of $e$, then any $P_4$ starting with~$e$ is rainbow.
If $u_1$ and $u_2$ are not the ends
of $e$, then there is a $P_4$ whose ends are $u_1$ and~$u_2$ and that 
$P_4$ must be rainbow. Either way, we obtain a contradiction. 

Since the graph is connected, it follows that every triangle is colored with 
the same three colors. Since this includes all the vertices, the result follows.
\end{Proof}

For example, it follows that if $G$ is a 
maximal outerplanar graph, then it follows that
$W^+(G; K_2,P_4)=3$.

\subsubsection{Cubic Graphs}

There are many cubic graphs with $W^+(G; K_2,P_4)=3$. These include, for example,
the claw-free cubic graphs (equivalently the ones where every vertex is in a triangle).
See Theorem~\ref{t:pathTriangle}.

For the largest value of the parameter, computer evidence suggests:

\begin{conjecture}
If $G$ is a connected cubic graph of order $n$, then 
$W^+(G; K_2,P_4)\le n/2+1$, with equality exactly when $G$ is bipartite.
\end{conjecture}

Certainly, by Observations \ref{o:bipartiteP4} and~\ref{o:matchingP4}
(and the fact that regular bipartite graphs have perfect matchings),
that value is obtained for all bipartite graphs.

\subsection{Forbidding rainbow $C_4$}

We conclude this section by considering proper colorings without rainbow
$4$-cycles. 

\begin{ob}
If $G$ is a maximal outerplanar graph, then $W^+(G; K_2, C_4 ) = 3$.
\end{ob}
\begin{Proof}
Consider two triangles sharing an edge. Then to avoid a rainbow $C_4$,
the two vertices not on the edge must have the same color. It follows that all
triangles have the same three colors.
\end{Proof}

In particular, we  again look at cubic graphs.
The parameter $W^-(G; K_2,C_4)$ for a cubic graphs $G$ is uninteresting:
the $3$-coloring provides such a WORM coloring, and so the parameter
is determined by whether $G$ is bipartite or not. Further, the
upper bound for $W^+(G; K_2,C_4)$ is trivial: one can have a cubic graph 
without a $4$-cycle.

Computer evidence suggests that:

\begin{conjecture}
If $G$ is a connected cubic graph of order $n$, then 
$W^+(G; K_2,C_4) \ge n/2$.
\end{conjecture}

This lower bound is achievable. Define a \emph{prism} as the cartesian product
of a cycle with~$K_2$. For $n$ even, a \emph{Mobius ladder} is defined by taking
the cycle on $n$ vertices and joining every pair of opposite vertices. Note that a prism
is bipartite when $n$ is a multiple of $4$, and a Mobius ladder is bipartite
when $n$ is not a multiple of $4$.

\begin{ob}
If $G$ is a nonbipartite Mobius ladder or prism of order $n$, 
then it holds that $W^+(G; K_2,C_4) = n/2$.
\end{ob}
\begin{Proof}
We first exhibit the coloring. Let $m=n/2$.
Say the vertices of the prism are $u_1, \ldots, u_m$ and $v_1, \ldots, v_m$, where
$u_i$ has neighbors $u_{i-1}$, $u_{i+1}$, and $v_i$ (arithmetic modulo~$m$) and similarly
for $v_i$. Then 
for $1\le i \le m$, give vertices $u_i$ and $v_{i+1}$ color $i$.

Say the vertices of the Mobius ladder are $w_1, \ldots, w_n$ where
$w_i$ has neighbors $w_{i-1}$, $w_{i+1}$, and $w_{i+m}$ (arithmetic modulo~$n$). Then 
for $2\le i \le m$, give vertices $w_i$ and $w_{i+m-1}$ color $i$, give vertex $w_1$
color $1$ and give vertex $w_{n-1}$ color $2$. For example, the coloring for the case $n=12$ 
is shown in Figure~\ref{f:mobiusColoring}.

\begin{figure}[!ht]
\centerline{\includegraphics{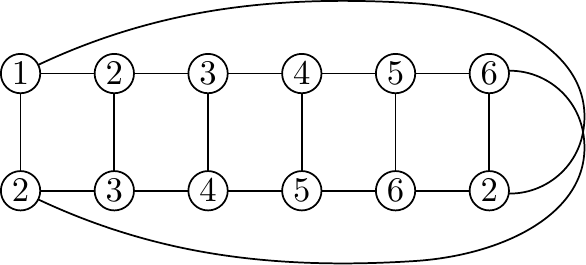}}
\caption{Coloring of Mobius ladder}
\label{f:mobiusColoring}
\end{figure}

Now, for the upper bound,
consider a color that is used $r$ times. A color bi-covers
a copy of $C_4$ if it contains vertices from consecutive rungs (where a rung is
an edge in two~$C_4$'s). Since the graph is not bipartite, the color cannot be
present in every rung. It follows that it can bi-cover at most $r-1$ copies of $C_4$.
Now, there are $m$ copies of $C_4$ (note that the prism of $C_4$ is bipartite so excluded).
It follows from Lemma~\ref{l:breakBound} that the number of colors is at most
$n-n/2=n/2$.
\end{Proof}

It appears that this extremal graph is unique for all orders.
 
\section{Other Results} \label{s:other}

\subsection{Paths and paths}

The natural strategy to color a long path without a rainbow $P_r$ yields
the following:

\begin{ob}
For any $m\ge 3$, it holds that
$W^+( P_n; P_m, P_r) = r^+( P_n; P_r ) = \lfloor (r-2)n/(r-1) \rfloor +1$.
\end{ob}
\begin{Proof}
Give the first $r-1$ vertices different colors, then the next vertex
the same color as the previous vertex, then the next $r-2$ vertices different
colors, and so on. This coloring has a monochromatic $P_2$ but not a monochromatic $P_3$,
and is easily seen to be best possible (as every copy of $P_r$ must
contain two vertices of the same color). 
\end{Proof}

\subsection{Bicliques and bicliques}

Next we revisit the case that $G$, $M$, and $R$ are bicliques. For $n\ge b$ it was proved that
$W^+(K_{n,n}; K_{1,b}, K_{1,b})=2b-2$ in~\cite{GWX15firstWorm} and that
$W^+(K_{n,n}; K_{b,b}, K_{b,b})=n+b-1$ in~\cite{GWX14secondWorm}. The case for stars is
special, but it is straight-forward to generalize the latter:

\begin{theorem}
Let $m\le n$ and $2\le a\le b$ with
$m\ge a$ and $n\ge b$. Then
\[
   r^+(K_{m,n}; K_{a,b}) = \max( a+n-1, \, b-1+\min(m,b-1) ) .
\]
\end{theorem}
\begin{Proof}
Consider a coloring $K_{m,n}$ without a rainbow $K_{a,b}$ and 
assume there are at least $a+b$ colors.
If one partite set has at least $a$ colors and the other partite
set has at least~$b$ colors, then one can choose $a$ colors from
the one and $b$ from the other that are disjoint and thus obtain
a rainbow $K_{a,b}$. So: either (1) there is a partite set
that has at most $a-1$ colors, or (2) both partite sets have
at most $b-1$ colors. In the first case, the maximum number
of colors possible is $a+n-1$. In the second case,
the maximum number of colors possible is $b-1+\min(m,b-1)$.
The theorem follows.
\end{Proof}

Note that in the above proof, the optimal number of colors can be
achieved by making the sets of colors in the two partite sets
disjoint. Thus, one obtains a similar value
for $W^+(K_{m,n}; M, K_{a,b})$ where $M$ is any nontrivial biclique.

\subsection{Grids without rainbow $4$-cycles}

We conclude this section with a result about forbidden $4$-cycles.
This result establishes a conjecture proposed in~\cite{GWX14secondWorm}.
Let $G_{m,n}$ denote the grid formed by the cartesian product of $P_m$ and $P_n$.

\begin{ob} \label{o:2times2grids}
For any grid and $s>0$,
$\BREAK{C_4}{s} \le 2(s-1)$.
\end{ob}
\begin{Proof}
We prove this bound by induction. Let $S$ be a set of $s$ vertices.
The bound is immediate when $S$ is contained in only one row. 
Now suppose $S$ intersects at least two rows. Let $S_1$ be a maximal set of consecutive
vertices of $S$ in the topmost row of $S$.
By the induction hypothesis,
the number of~$C_4$'s that contain at least two vertices in $S \setminus S_1$ is at
most $2(|S|-|S_1|)-2$. Further, the number of $C_4$'s that contain at least one vertex in~$S_1$ 
and least two vertices in $S$ is at most $2|S_1|$: there are $|S_1|-1$ possible copies above $S_1$
and at most $|S_1|+1$ copies below. Hence, the number of $C_4$'s 
that $S$ bi-covers is at most $2(|S|-|S_1|)-2 + 2|S_1| = 2|S|-2$.
\end{Proof}

In~\cite{GWX14secondWorm} a $(C_4,C_4)$-WORM coloring is given and it is conjectured
that this is best possible. This we now show:

\begin{theorem} \label{t:C4inGrids}
For $G_{m,n}$ the $m\times n$ grid, it holds that $W^+(G_{m,n}; C_4, C_4 ) = \lfloor (m+1)(n+1)/2 \rfloor - 1$.
\end{theorem}
\begin{Proof}
The lower bound was proved in~\cite{GWX14secondWorm}. 
The upper bound follows from Lemma~\ref{l:breakBound} and Observation~\ref{o:2times2grids}:
There are $(m-1)(n-1)$ copies of $C_4$, and so $r^+(G_{m,n}; C_4) 
\le mn - (m-1)(n-1)/2  = (m+1)(n+1)/2-1$.
\end{Proof}

\section{Other Directions}

We conclude with some thoughts on future directions. Apart from the specific
open problems raised here, a direction that looks interesting is the case where $M$ and
$R$ are both stars. Also of interest is where the host graph is a product graph.

\bibliographystyle{amsplain}
\bibliography{genWorm}

\end{document}